\input amstex
\input amsppt.sty   
\hsize 30pc
\vsize 47pc
\def\nmb#1#2{#2}         
\def\cit#1#2{\ifx#1!\cite{#2}\else#2\fi} 
\def\totoc{}             
\def\idx{}               
\def\ign#1{}             

\redefine\o{\circ}
\define\X{\frak X}
\define\al{\alpha}

\define\ga{\gamma}
\define\de{\delta}

\define\ze{\zeta}
\define\et{\eta}
\define\th{\theta}

\define\ka{\kappa}
\define\la{\lambda}
\define\rh{\rho}
\define\si{\sigma}

\define\ph{\varphi}

\define\ps{\psi}
\define\om{\omega}
\define\Ga{\Gamma}
\define\De{\Delta}
\define\Th{\Theta}
\define\La{\Lambda}
\define\Si{\Sigma}
\define\Ph{\Phi}
\define\Ps{\Psi}
\define\Om{\Omega}
\redefine\i{^{-1}}
\define\x{\times}
\define\row#1#2#3{#1_{#2},\ldots,#1_{#3}}
\define\End{\operatorname{End}}
\define\Fl{\operatorname{Fl}}

\define\Ad{\operatorname{Ad}}
\define\ad{\operatorname{ad}}
\define\pr{{\operatorname{pr}}}
\define\tr{{\operatorname{tr}}}
\redefine\L{{\Cal L}}
\define\ddt{\left.\tfrac \partial{\partial t}\right\vert_0}
\define\g{{\frak g}}
\define\h{{\frak h}}

\define\Conj{\operatorname{Conj}}
\def\today{\ifcase\month\or
 January\or February\or March\or April\or May\or June\or
 July\or August\or September\or October\or November\or December\fi
 \space\number\day, \number\year}
\topmatter
\title   Poisson structures on the cotangent bundle of a Lie group or 
a principle bundle and their reductions
\endtitle
\author  D\. Alekseevsky\\
J\. Grabowski\\
G\. Marmo\\
P\. W\. Michor  \endauthor
\affil
\leftheadtext{\smc Alekseevsky, Grabowski, Marmo, Michor}
\rightheadtext{\smc Poisson structures}
Erwin Schr\"odinger International Institute of Mathematical Physics, 
Wien, Austria
\endaffil
\address
Erwin Schr\"odinger International Institute of Mathematical Physics, 
Pasteurgasse 6/7, A-1090 Wien, Austria
Wien, Austria
\endaddress
\address 	 D\. V\. Alekseevsky: 
gen. Antonova 2 kv 99, 117279 Moscow B-279, Russia
\endaddress
\address J\. Grabowski: Institute of Mathematics, 
University of Warsaw,\newline 
ul. Banacha 2, PL 02-097 Warsaw, Poland
\endaddress
\email jagrab\@mimuw.edu.pl \endemail
\address
G. Marmo:
Dipart\. di Scienze Fisiche - Universit\`a di Napoli,
Mostra d'Oltremare, Pad.19, I-80125 Napoli, Italy.
\endaddress
\email gimarmo\@na.infn.it \endemail
\address
P\. W\. Michor: Institut f\"ur Mathematik, Universit\"at Wien,
Strudlhofgasse 4, A-1090 Wien, Austria; and 
Erwin Schr\"odinger International Institute of Mathematical Physics, 
Pasteurgasse 6/7, A-1090 Wien, Austria
\endaddress
\email michor\@pap.univie.ac.at, michor\@esi.ac.at \endemail
\date {\today} \enddate
\endtopmatter

\document

\heading Table of contents \endheading
\noindent Introduction \leaders \hbox to 1em{\hss .\hss }\hfill {\eightrm 1}\par 
\noindent 1. Liouville 1-forms on fiber bundles and the lifting of vector fields \leaders \hbox to 1em{\hss .\hss }\hfill {\eightrm 3}\par 
\noindent 2. The canonical Poisson structure on $T^*G$ \leaders \hbox to 1em{\hss .\hss }\hfill {\eightrm 5}\par 
\noindent 3. Generalizing momentum mappings \leaders \hbox to 1em{\hss .\hss }\hfill {\eightrm 11}\par 
\noindent 4. Symplectic structures on cotangent bundles of principal bundles \leaders \hbox to 1em{\hss .\hss }\hfill {\eightrm 13}\par 

\head\totoc Introduction \endhead
The standard description of physical systems,  both for particles and 
for fields, usually starts with an action principle. When particles 
are thought of as test particles one considers fields as given, 
i\.e\. as external fields, and only the point particle dynamics is 
considered. In this framework the Lagrangian function usually is the 
sum of three terms: A kinematic term which is quadratic in the 
velocities, a current-potential coupling term which is linear in the 
velocities, and a term which depends only on the positions like 
electrostatic potential. 

When one passes to the Hamiltonian description in the symplectic or 
Poisson formalism one may think of the magnetic field absorbed in a 
change of coordinates, so $p$ is replaced by $p+eA$, and the 
Hamiltonian takes into account the electrostatic potential or other 
effective potentials.

There are situations, however, like in the electric monopole system 
for instance, where the magnetic field cannot be absorbed in a change 
of coordinates. This is due to the fact that the associated 
symplectic structures are in different cohomology classes. From this 
point of view the Lagrangian formalism seems to be able  to take into 
account a specific identification of physical variables which is not 
possible in the other descriptions. Nevertheless Poisson brackets 
seem to be an important starting point for various quantization 
procedures. More recently a special class of Poisson brackets has 
been considered to represent the classical limits of quantum 
groups. It seems therefore appropriate to look at Poisson brackets in a 
more direct way in order to learn how one may be able to `add 
interactions' directly to the brackets. To this aim we would like to 
analyze the brackets that are appropriate to describe point 
particles carrying internal degrees of freedom, isospin-like 
variables, interacting with external Yang-Mills fields. Eventually 
one might be able to either deform these brackets or to learn how a 
given bracket can be read off as arising from reduction of a system 
described by an action principle. The present paper tries to 
provide partial answers to our questions while referring to future 
work for further developments.

The paper is organized as follows. In section \nmb!{1} we define a 
Liouville form on a symplectic manifold $(E,\om)$ fibered over a 
manifold $M$ with Lagrangian fibers as a horizontal 1-form $\Th$ with 
$d\Th=\om$. We denote by $\frak A(\Th)$ the Lie algebra of all 
projectable vector fields on $E$ which preserve the Liouville form 
$\Th$. A right inverse to the   projection $\frak A(\Th)\to \X(M)$ is 
called the Liouville lift $\X(M)\to \frak A(\Th)$. It assigns to a 
vector field $X\in\X(M)$ the Hamiltonian vector field $X^\Th$ on $E$ 
for the function $-\Th(X\o p)$, where $p:E\to M$ is the projection. 
Using the Liouville lift we obtain the explicit expression for the 
Poisson structure $\La=\om\i$. We also sketch a reverse of this 
construction.

In section \nmb!{2} we apply these simple construction to the 
cotangent bundle $\pi:T^*G\to G$ of a Lie group $G$ and we describe 
explicitly the standard symplectic form $\om$ on $T^*G$ as 
follows:
$$
\om=d\Th = \frac12\left(\langle d\ze^l,\pi^*\ka^r 
     \rangle^\wedge_\g + \langle d\ze^r,\pi^*\ka^l 
     \rangle^\wedge_\g \right),
$$
where $\ka^l$, $\ka^r$ are the left and right Maurer-Cartan forms on 
$G$, and where $\ze^l$, $\ze^r$ are the momenta of the left and right 
action of $G$ onto itself. The associated Poisson structure is given by
$$
\La = \frac12\sum_i \left( R^*_i\wedge Z^r_i + L^*_i\wedge Z^l_i 
\right),
$$
where $L_i^*$ and $R_i^*$ are the flow lifts to $T^*G$ of the left 
and right invariant vector fields on $G$ corresponding to the basis 
vectors $X_i$ of the Lie algebra $\g$, and where the vertical vector 
fields $Z_j^l$, $Z_j^r$ are defined to be $\om$-dual to the 1-forms 
$-\pi^*\ka^l_j$, $-\pi^*\ka^r_j$ which are pullbacks of the components 
of the Maurer-Cartan forms: $\ka^l=\sum\ka^l_j\otimes X_j$ and 
$\ka^r=\sum\ka^r_j\otimes X_j$.

Some generalization of the above construction is presented in section 
\nmb!{3}. The starting point is that instead of the standard momentum 
mapping we consider an arbitrary smooth mapping $f:T^*G\to \g^*$ to 
define a generalized Liouville form 
$$\Th_f = \langle f,\pi^*\ka^l\rangle \in \Om^1(T^*G)$$
and a derived closed 2-form $\om_f=d\Th_f$. In general the 2-form 
$\om_f$ is degenerate and we need a reduction in order to obtain a 
symplectic form. We consider some examples. In particular, for a 
compact Lie group $G$, starting with a reduced momentum mapping we 
come after reduction to a symplectic form and associated Poisson 
structure $\La$ on the manifold $G\x C$, where $C$ is an open 
Weyl-chamber in the dual to a Cartan subalgebra. This Poisson 
structure was also described in \cit!{1}.

In the last section we generalize the construction of section 
\nmb!{3} to the case of a principal bundle $p:P\to M$ with structure 
group $G$. We choose a principal connection with connection form 
$\ga:TP\to \g=\operatorname{Lie}(G)$ and note that its pullback to 
the cotangent bundle $\pi_P:T^*P\to P$ gives a $\g$-valued 
$G$-equivariant $\pi_P$-horizontal 1-form $\pi^*\ga:T(T^*P)\to \g$. 
We define a generalized momentum mapping as an arbitrary 
$G$-equivariant $\g^*$-valued function $f:T^*P\to \g^*$ and consider 
the generalized Liouville form
$$
\Th_f = \langle f,\pi^*\ga\rangle = \sum f_j.\ga^j\in \Om^1(T^*P).
$$
The corresponding 2-form $\om_f=d\Th_f$ on $T^*P$ is degenerate in 
general. However, we can try to reduce it to a symplectic form on the 
reduced manifold $P/\ker\om_f$. As a simple example we consider the 
form $\Th_f$ associated to the standard momentum mapping 
$f:T^*P\to \g^*$ associated to the Hamiltonian right action of $G$ on 
$T^*P$.
We call it 
the vertical Liouville form. We also calculate $d\Th_f$ in 
coordinates. Using local trivializations of the bundle $p:P\to M$ we 
consider a Liouville form $\Th_\ga$ which is the sum of the vertical 
one and a horizontal one which is the pullback of the standard one on 
$T^*M$. The associated 2-form $\om_\ga=d\Th_\ga$ is non degenerate 
and we calculate the inverse Poisson structure $\La_\ga$ in 
coordinates. Since $\La_\ga$ is $G$-invariant we can factorize it and 
as a result we obtain a Poisson structure $\tilde\La_\ga$ on the 
orbit space $T^*P/G$, which is degenerate and does not come from a 
symplectic  form.  The Poisson bracket associated to $\tilde\La_\ga$ 
was considered in \cit!{9}.

Then we consider the case of a trivial principal bundle 
$P=\Bbb R^n\x G$ with a compact structure group $G$ and take as a 
momentum mapping $f$ the projection of the canonical momentum 
mapping onto the dual $\h^*$ of a Cartan subalgebra $\h$ of 
$\g=\operatorname{Lie}(G)$. We show that in this case our 
construction leads to a Poisson structure on 
$T^*\Bbb R^n\x G\x \h^*$ which has singularities on the walls of the 
Weyl chambers in $\h^*$.

We consider also the case when a principal bundle $p:P\to M$ is 
equipped with a displacement form (soldering form) 
$\th:TP\to \Bbb R^n$ and describes a $G$-structure on $M$. In this 
case any $G$-equivariant function $f:T^*P\to \Bbb R^n\x\g^*$ defines 
a $G$-invariant generalized Liouville form on $T^*P$.

\head\totoc\nmb0{1}. Liouville 1-forms on fiber bundles and the 
lifting of vector fields \endhead

\subhead\nmb.{1.1}. Liouville forms on fiber bundles \endsubhead
Let $p:E\to M$ be a locally trivial smooth fiber bundle, and let 
$\om\in\Om^2(E)$ be a symplectic form on $E$. A 1-form 
$\Th\in\Om^1(E)$ is called a \idx{\it Liouville form} if: 
\roster
\item It is horizontal: $i_Y\Th=0$ for $Y$ in the vertical bundle 
       $VE:=\ker(Tp:TE\to TM)$. 
\item $\om = d\Th$.
\item The fibers of $p:E\to M$ are Lagrangian submanifolds.
\endroster

\proclaim{Lemma} In this situation, if a Liouville form $\Th$ exists, 
$(E,\Th)$ is locally fiber respecting diffeomorphic to 
$(T^*M,\Th_M)$, where $\Th_M$ is the canonical Liouville form on 
$T^*M$. 
\endproclaim

\demo{Proof}
For a point $u\in E$, choose local coordinates $q^i$ near 
$p(u)\in M$. Then since $\Th$ is horizontal, 
near $u$ we have $\Th = \sum_ip_i\,dq^i$ for local smooth 
functions $p_i$ on $E$. Since $\om=d\Th=\sum_idp_i\wedge dq^i$ is 
symplectic, and since the fibers of $p:E\to M$ are Lagrangian 
submanifolds, $dq^1,\dots,dq^n,dp_1,\dots,dp_n$ is a local frame for 
$T^*E$, so $q^1,\dots,q^n,p_1,\dots,p_n$ is a coordinate system near 
$u$ on $E$.
\qed\enddemo

\subhead\nmb.{1.2}. Liouville lift of a vector field  \endsubhead
Let $p:E\to M$ be a fiber bundle with a symplectic form 
$\om\in\Om^2(E)$ and a Liouville form $\Th\in\Om^1(E)$.
Since $\Th: TE\to \Bbb R$ is horizontal it factors to a form on the 
quotient bundle $TE/VE\cong p^*TM\ \to \Bbb R$, and for each vector 
field $X\in\X(M)$ on $M$ we may consider $\Th(X\o p)$ as a function on 
$E$. 

A vector field $X^\Th\in\X(E)$ is called the \idx{\it Liouville lift} 
of a vector field $X\in\X(M)$ if $X^\Th$ is  
projectable onto $X$ and preserves $\Th$, that is $Tp\o X^\Th= X\o p$ 
and $\L_{X^\Th}\Th=0$. 
But then 
$$0 = \L_{X^\Th}\Th = i_{X^\Th}\om + di_{X^\Th}\Th = 
     i_{X^\Th}\om + d(\Th(X\o p)).$$
Hence $X^\Th$ is the Hamiltonian vector field $-H_{\Th(X\o p)}$
for the function $-\Th(X\o p)\in C^\infty(E,\Bbb R)$.

Using the fact from lemma \nmb!{1.1} that $(E,\Th)$ is locally 
diffeomorphic to $(T^*M,\Th_M)$ one can check easily that 
$X^\Th=H_{-\Th(X\o p)}$ is indeed projectable onto $X\in \X(M)$. It 
is also well known for cotangent bundles that the mapping 
$\X(M)\to C^\infty(T^*M,\Bbb R)$ given by 
$X\mapsto -\Th(X)$ is an injective 
homomorphism of Lie algebras from the Lie 
algebra of vector fields to the Lie algebra of functions with the 
standard Poisson bracket, so this holds also 
for the general situation.

\subhead\nmb.{1.3}. A standard frame of $E$ \endsubhead
Let us assume now that the manifold $M$ is parallelizable and that 
$X_1,\dots,X_n$ form a global frame field for $TM$. Let us denote by 
$\al_1,\dots,\al_n\in\Om^1(M)$ the dual coframe field.
Then we consider the frame field 
$$X_1^\Th,\dots,X_n^\Th, Z_1,\dots Z_n,$$
where $X_i^\Th$ is the Liouville lift from \nmb!{1.2} and where 
$$Z_j:=\om\i p^*\al_j,\quad\text{ or  equivalently }i_{Z_j}\om=p^*\al_j \tag1$$
is the vertical vector field dual to the 1-form 
$p^*\al_i\in\Om^1(E)$. From lemma \nmb!{1.1} and a computation in 
$T^*M$ it follows that $[Z_i,Z_j]=0$.

Then the Liouville form $\Th$ and the symplectic form $\om=d\Th$ 
may be written as 
$$\align
\Th &= \sum_i f_i p^*\al_i\quad\text{ for }f_i\in C^\infty(E,\Bbb R)\\
\om &= d\Th = \sum_i df_i\wedge p^*\al_i + \sum_i f_i\,p^*d\al_i \tag2\\
&= \sum_i df_i\wedge p^*\al_i + \sum_{i,j,k} 
     f_ic^i_{jk}p^*\al_j\wedge p^*\al_k,
\endalign$$
where the \idx{\it torsion functions} 
$c^i_{jk}\in C^\infty(M,\Bbb R)$ are given by 
$[X_j,X_k]=\sum_i c^i_{jk}X_j$. 
 From the definition of $X_i^\Th$ we have 
$$
i_{X_j^\Th}\om = - d\Th(X_j\o p) = - df_j,
$$
thus the Poisson structure associated to $\om$ is given by 
$$
\La=\om\i=\sum_i X_j^\Th\wedge Z_i + \sum_{i,j,k} 
     f_ic^i_{jk}Z_j\wedge Z_k \in C^\infty(\La^2TE)\tag 2
$$

\subhead\nmb.{1.4}. The reversed construction \endsubhead
Let $p:E\to M$ be a fiber bundle whose fiber dimension equals the 
dimension of the base space $M$. Let $M$ be parallelizable and let 
$$\al = (\al_1,\dots,\al_n): TM \to \Bbb R^n =: V$$
is a coframe field for $M$. If we are given a $V^*$-valued function 
$f:E\to V^*$, we may consider the 1-form 
$$\Th_f:= \langle f, p^*\al\rangle\in \Om^1(E).$$

\proclaim{Proposition} The 1-form $\Th_f$ defines a symplectic 
structure $\om_f=d\Th_f$ on $E$ if and only if $f|E_x:E_x\to V^*$ is 
a local diffeomorphism for each fiber $E_x$, $x\in M$. 
In this case, $\Th_f$ is a Liouville form for $(E, \om_f)$. The 
associated Poisson bivector field $\La_f=\om_f\i$ is given by formula 
\nmb!{1.3}.\thetag2, \qed 
\endproclaim

\head\totoc\nmb0{2}. The canonical Poisson structure on $T^*G$ \endhead

\subheading{\nmb.{2.1} Products of differential forms}
Let $\rh:\g\to \frak g\frak l(V)$ be a representation of a Lie 
algebra $\g$ in a finite dimensional 
vector space $V$ and let $M$ be a smooth manifold.

For vector valued differential forms 
$\ph\in\Om^p(M;\g)$ and $\Ps\in\Om^q(M;V)$ we define the 
form	$\rh^\wedge (\ph)\Ps\in\Om^{p+q}(M;V)$	by 
$$\multline
(\rh^\wedge (\ph)\Ps)(\row X1{p+q}) =\\
= \frac 1{p!\,q!} \sum_{\si} \text{sign}(\si)
 	\rh(\ph(\row X{\si1}{\si p}))\Ps(\row X{\si(p+1)}{\si(p+q)}).
\endmultline$$
Then $\rh^\wedge (\ph):\Om^*(M;V)\to \Om^{*+p}(M;V)$ is a graded 
$\Om(M)$-module homomorphism of degree $p$.

Recall also that $\Om(M;\g)$ is a graded Lie algebra with the bracket 
$[\quad,\quad]^\wedge = [\quad,\quad]^\wedge_\g$ given by
$$\multline [\ph,\ps]^\wedge(\row X1{p+q}) = \\
	= \frac 1{p!\,q!} \sum_{\si} \text{sign}\si\,
	[\ph(\row	X{\si1}{\si p}),\ps(\row X{\si(p+1)}{\si(p+q)})]_{\g},
\endmultline $$
where $[\quad,\quad]_{\g}$ is the bracket in $\g$.
One may easily check that for the graded commutator in 
$\End(\Om(M;V))$ we have
$$\rh^\wedge ([\ph,\ps]^\wedge ) = 
[\rh^\wedge (\ph),\rh^\wedge (\ps)] = 
\rh^\wedge (\ph)\o \rh^\wedge (\ps) - (-1)^{pq} 
\rh^\wedge (\ps)\o \rh^\wedge (\ph)$$
so that $\rh^\wedge :\Om^*(M;\g) \to \End^*(\Om(M;V))$ is a 
homomorphism of graded Lie algebras.

For any vector space $V$ let
$\bigotimes V$ be the tensor algebra generated by $V$.
For $\Ph,\Ps\in \Om(M;\bigotimes V)$ we will use the associative 
bigraded product
$$\multline
(\Ph\otimes_\wedge \Ps)(\row X1{p+q}) = \\
= \frac 1{p!\,q!} \sum_{\si} \text{sign}(\si)
 	\Ph(\row X{\si1}{\si p})\otimes \Ps(\row X{\si(p+1)}{\si(p+q)})
\endmultline$$

In the same spirit we will use the following product:
Let $V$ be a finite dimensional vector space with dual $V^*$, and let 
$\langle  \quad,\quad\rangle:V^*\x V\to \Bbb R$ be the duality 
pairing. 
For $\Ph\in \Om^p(M;V^*)$ and $\Ps\in \Om^q(M;V)$ we consider the 
`product' $\langle \Ph,\Ps\rangle^\wedge \in \Om(M)$ which is given 
by
$$\multline
\langle \Ph,\Ps\rangle^\wedge (\row X1{p+q}) = \\
= \frac 1{p!\,q!} \sum_{\si} \text{sign}(\si)
 	\langle \Ph(\row X{\si1}{\si p}), 
     \Ps(\row X{\si(p+1)}{\si(p+q)}\rangle.
\endmultline$$

\subhead\nmb.{2.2}. Notation for Lie groups
\endsubhead
Let $G$ be a Lie group with Lie algebra $\g=T_eG$, 
multiplication $\mu:G\x G\to G$, and for $g\in G$ 
let $\mu_g, \mu^g:G\to G$ denote the left and right translation, 
$\mu(g,h)=g.h=\mu_g(h)=\mu^h(g)$. 

Let $L,R:\g\to \X(G)$ be the left 
and right invariant vector field mappings, given by 
$L_X(g)=T_e(\mu_g).X$ and $R_X=T_e(\mu^g).X$, respectively. 
They are related by $L_X(g)=R_{\Ad(g)X}(g)$.
Their flows are given by 
$$\Fl^{L_X}_t(g)= g.\exp(tX)=\mu^{\exp(tX)}(g),\quad
\Fl^{R_X}_t(g)= \exp(tX).g=\mu_{\exp(tX)}(g).$$

Let $\ka^l,\ka^r:\in\Om^1(G,\g)$ be the left and right Maurer-Cartan 
forms, given by $\ka^l_g(\xi)=T_g(\mu_{g\i}).\xi$ and 
$\ka^r_g(\xi)=T_g(\mu^{g\i}).\xi$, respectively. These are the 
inverses to $L,R$ in the following sense: $L_g\i=\ka^l_g:T_gG\to\g$ 
and $R_g\i=\ka^r_g:T_gG\to\g$. They are related by 
$\ka^r_g=\Ad(g)\ka^l_g:T_gG\to\g$ and they satisfy the Maurer-Cartan 
equations $d\ka^l+\frac12[\ka^l,\ka^l]^\wedge =0$ and 
$d\ka^r-\frac12[\ka^r,\ka^r]^\wedge =0$. 

The (exterior) derivative of the function $\Ad:G\to GL(\g)$ can be 
expressed by
$$d\Ad = \Ad.(\ad\o\ka^l) = (\ad\o \ka^r).\Ad,$$
which follows from
$d\Ad(T\mu_g.X) = \frac d{dt}|_0 \Ad(g.\exp(tX))
= \Ad(g).\ad(\ka^l(T\mu_g.X))$.

We also consider the left moment mapping $\ze^l:T^*G\to\g^*$ which is 
related to the left action of $G$ on $T^*G$ by left translations, 
given by $\ze^l_g = ((\ka^r_g)\i)^*$. Note here the transition from 
left to right: right invariant vector fields generate left 
translations.
Similarly we consider the right moment mapping $\ze^r:T^*G\to\g^*$ which is 
related to the right action of $G$ on $T^*G$ by right translations, 
given by $\ze^r_g = ((\ka^l_g)\i)^*$. They are related by 
$\ze^l_g=\Ad(g\i)^*\ze^r_g=\Ad^*(g).\ze^r_g:T^*_gG\to\g^*$.

\subhead\nmb.{2.3}. The canonical symplectic structure on $T^*G$  
\endsubhead
We consider now the tangent bundle $\pi_G:TG\to G$ and the  
cotangent bundle $\pi=\pi_G:T^*G\to G$ of the Lie group $G$. We use 
$\pi$ for the projection of each bundle which is derived from the 
tangent bundle in a direct way.  On $T^*G$ we consider the Liouville 
form $\Th:T(T^*G)\to \Bbb R$, $\Th\in\Om^1(T^*G)$, which is given by 
$$
\Th(\Xi) = \langle \pi_{T^*G}(\Xi), T(\pi_G).\Xi\rangle_{TG},\tag1
$$
where we use the following commutative diagram
$$\CD
TT^*G @>{T(\pi_G)}>> TG \\
@V{\pi_{T^*G}}VV     @V{\pi_G}VV\\
T^*G  @>>{\pi}>         G.
\endCD$$
Considering the momentum mapping $\ze^l:T^*G\to\g^*$ as a function and 
the pullback $\pi_{T^*G}^*\ka^r=\pi^*\ka^r\in\Om^1(T^*G,\g)$ 
of the right Maurer-Cartan 
form, we have for 
$\Xi\in TT^*G$ with `lowest footpoint' $g\in G$:
$$\align
\langle \ze^l,\pi^*\ka^r\rangle_\g(\Xi) &= 
     \langle ((\ka^r_g)\i)^*(\pi_{T^*G}(\Xi)),\ka^r_g(T(\pi_G).\Xi) \rangle_\g \\
& =\langle \pi_{T^*G}(\Xi),(\ka^r_g)\i\ka^r_g(T(\pi_G).\Xi) \rangle_\g \\
&= \langle \pi_{T^*G}(\Xi),T(\pi_G).\Xi \rangle_\g =\Th(\Xi).
\endalign$$
Similarly we have $\Th = \langle \ze^r,\pi^*\ka^l \rangle_\g $ 
and thus also 
$$\Th= \tfrac12\left(\langle \ze^l,\pi^*\ka^r \rangle_\g + 
\langle \ze^r,\pi^*\ka^l \rangle_\g\right).\tag2$$
Let us now compute the exterior derivative of the one form 
$\Th=\langle \ze^l,\pi^*\ka^r \rangle_\g$ in order to get an 
expression for the symplectic structure $\om\in\Om^2(T^*G)$:
$$\align
\om&=d\Th = d(\langle \ze^l,\pi^*\ka^r \rangle_\g) 
= \langle d\ze^l,\pi^*\ka^r \rangle^\wedge_\g 
     + \langle \ze^l,\pi^* d\ka^r \rangle_\g\\
&= \langle d\ze^l,\pi^*\ka^r \rangle^\wedge_\g 
     -\tfrac12 \langle \ze^l,\pi^* [\ka^r,\ka^r]^\wedge_\g  
     \rangle_\g.
\endalign$$
by the Maurer-Cartan equation. Similarly we get 
$$ \om=d\Th = d(\langle \ze^r,\pi^*\ka^l \rangle_\g) 
     = \langle d\ze^r,\pi^*\ka^l \rangle^\wedge_\g 
     +\tfrac12 \langle \ze^r,\pi^* [\ka^l,\ka^l]^\wedge_\g  
     \rangle_\g. $$
For $\Xi_i\in T_\xi(T^*G)$ with $\pi(\xi)=g\in G$ we have
$$\align
\langle \ze^r,\pi^* [\ka^l,\ka^l]^\wedge_\g\rangle_\g(\Xi_1,\Xi_2) 
&= \langle \ze^r(\xi),[\pi^* \ka^l(\Xi_1),
     \pi^* \ka^l(\Xi_2)]_\g\rangle_\g\\
&= \langle\Ad(g\i)^* \ze^r(\xi),\Ad(g)[\pi^* \ka^l(\Xi_1),
     \pi^* \ka^l(\Xi_2)]_\g\rangle_\g\\
&= \langle\Ad^*(g) \ze^r(\xi),[\Ad(g).\ka^l(T\pi.\Xi_1),
     \Ad(g).\ka^l(T\pi.\Xi_2)]_\g\rangle_\g\\
&= \langle \ze^l(\xi),[\ka^r(T\pi.\Xi_1),
     \ka^r(T\pi.\Xi_2)]_\g\rangle_\g\\
&= \langle \ze^l,\pi^* 
     [\ka^r,\ka^r]_\g\rangle_\g(\Xi_1,\Xi_2). 
\endalign$$
 From \thetag2 we get the following final formula for the canonical 
symplectic structure $\om$ on $T^*G$:
$$
\om=d\Th = \frac12\left(\langle d\ze^l,\pi^*\ka^r 
     \rangle^\wedge_\g + \langle d\ze^r,\pi^*\ka^l 
     \rangle^\wedge_\g \right).\tag3
$$

\subhead\nmb.{2.4}. The canonical symplectic structure on $T^*G$  in 
coordinates \endsubhead
We fix now a basis $X_1,\dots, X_n$ of the Lie algebra $\g$, with 
dual basis $\xi_1,\dots,\xi_n$ in $\g^*$, and with structure constants 
$[X_i,X_j]=\sum_kc^k_{ij}X_k$. We may expand the Maurer-Cartan forms 
and the moment mappings in terms of these bases, as
$$\alignat2
\ka^l &= \sum_i \ka^l_i.X_i,&\quad 
     \ka^l_i:&= \langle \xi_i,\ka^l\rangle \in \Om^1(G),\\
\ka^r &= \sum_i \ka^r_i.X_i,&\quad 
     \ka^r_i:&= \langle \xi_i,\ka^r\rangle \in \Om^1(G),\\
\ze^l &= \sum_i \ze^l_i.\xi_i,&\quad 
     \ze^l_i:&= \langle \ze^l,X_i\rangle \in C^\infty(T^*G,\Bbb R),\\
\ze^r &= \sum_i \ze^r_i.\xi_i,&\quad 
     \ze^r_i:&= \langle \ze^r,X_i\rangle \in C^\infty(T^*G,\Bbb R).
\endalignat$$
The Maurer-Cartan equations then become
$$\gather 
d\ka^l_i + \tfrac12\sum_{jk}c^i_{jk} \ka^l_j\wedge \ka^l_k = 0,\\
d\ka^r_i - \tfrac12\sum_{jk}c^i_{jk} \ka^r_j\wedge \ka^r_k = 0.
\endgather$$
The Liouville form then becomes 
$$\Th = \sum_i \ze^l_i.\pi^*\ka^r_i = \sum_i \ze^r_i.\pi^*\ka^l_i 
     = \tfrac12\sum_i \left(\ze^l_i.\pi^*\ka^r_i 
     + \ze^r_i.\pi^*\ka^l_i \right) $$
and the symplectic form looks like
$$\align
\om &= \sum_i \left( d\ze^l_i\wedge \pi^*\ka^r_i 
     + \frac12 \ze^l_i\sum_{jk} c^i_{jk} \pi^*\ka^r_j\wedge \pi^*\ka^r_k \right)\\
&= \sum_i \left( d\ze^r_i\wedge \pi^*\ka^l_i 
     - \frac12 \ze^r_i\sum_{jk} c^i_{jk} \pi^*\ka^l_j\wedge \pi^*\ka^l_k \right)\\
&= \frac12\sum_i \left( d\ze^l_i\wedge \pi^*\ka^r_i 
     + d\ze^r_i\wedge \pi^*\ka^l_i \right).
\endalign$$
Next we consider for the basis vectors $X_i$ of $\g$ the left 
invariant vector fields $L_i:=L_{X_i}$ and the right invariant vector 
fields $R_i=R_{X_i}$ on $G$ and their `flow lifts' to $T^*G$ which 
are given by $L^*_i := \ddt T^*(\Fl^{L_i}_t)$  and by 
$R^*_i := \ddt T^*(\Fl^{R_i}_t)$. 
Then $L^*_i\in\X(T^*G)$ is $\pi$-related to $L_i\in\X(G)$, similarly 
$T(\pi)\o R^*_i = R_i\o \pi: T^*G\to TG$. Thus 
$$\align
\pi^*\ka^l_i(L_j^*)&= \ka^l_i(T(\pi).L^*_j)=\ka^l_i(L_j)\o\pi = 
     \de_{ij},\\
\pi^*\ka^r_i(R_j^*)&= \de_{ij}.
\endalign$$
By general principles 
their flows preserve the Liouville form $\Th$,
$$0=\L_{L^*_j}\Th = i_{L^*_j}d\Th + di_{L^*_j}\Th 
     = i_{L^*_j}d\Th + d\ze^r_j,$$
so we conclude that $L^*_j$ is the hamiltonian vector field 
for the generating function $\ze^r_j$.
Similarly  $-i_{R^*_j}\om=d\ze^l_j$. 

We consider the vector fields $Z^l_j, Z^r_j\in\X(T^*G)$ which are 
given by 
$$-i_{Z^l_j}\om = \pi^*\ka^l_j,\quad -i_{Z^r_j}\om = \pi^*\ka^r_j.$$
The fields $Z^l_j$, $Z^r_j$ are vertical, i\.e\. in the kernel of 
$T(\pi)$.
Then the Poisson structure $\La=\om\i$ is given by 
$$\La = \frac12\sum_i \left( R^*_i\wedge Z^r_i + L^*_i\wedge Z^l_i \right)$$

\subheading{\nmb.{2.5}. The tangent group of a Lie group} 
As motivation for the following we recall that for a Lie group $G$ 
the tangent group $TG$
is also a Lie group with multiplication $T\mu$ and inversion $T\nu$, the 
tangent mapping of the inversion $\nu$ on $G$,
given by $T_{(a,b)}\mu.(\xi_a,\et_b) =
T_a(\mu^b).\xi_a + T_b(\mu_a).\et_b$ and $T_a\nu.\xi_a = -
T_e(\mu_{a\i}).T_a(\mu^{a\i}).\xi_a$.

\proclaim{Lemma} In the right trivialization, i\.e\. via the isomomorphism 
$(\ka^r,\pi):TG\to \g\x G$, the group structure on $TG$
looks as follows:
$$(X,a).(Y,b) = (X + \Ad(a)Y,a.b),\quad (X,a)\i = (-\Ad(a\i)X,a\i).$$
In the left trivialization, i\.e\. via the isomomorphism 
$(\pi,\ka^l):TG\to G\x \g$, the group structure on $TG$
looks as follows:
$$(a,X).(b,Y)=(ab,\Ad(b\i)X+Y),\quad (a,X)\i=(a\i,-\Ad(a)X).\qed$$
\endproclaim

\subhead\nmb.{2.6}. The Lie group $T^*G$ as semidirect product 
\endsubhead
The Lie group $G$ acts on the dual of its Lie algebra by the 
coadjoint representation $\Ad^*(g) = \Ad(g\i)^*$. So we can consider 
on $\g^*\x G$ the structure of a Lie group given by the semidirect 
product
$$(a,\xi).(b,\et)=(ab,\Ad^*(b\i)\xi+\et),\quad 
     (a,\xi)\i=(a\i,-\Ad^*(a)\xi).$$
We will consider on $T^*G$ the Lie group structure induced from this 
semidirect product by the diffeomorphism $(\ze^l,\pi):T^*G\to \g^*\x G$, 
whose inverse is given by $(\ze^l,\pi)\i(\al,g)= 
(T_g(\mu^{g\i}))^*\al = T^*(\mu^g)\al$.
This Lie group structure on $T^*G$ has no obvious intrinsic meaning, 
but it has nice relations to some structures which are naturally 
given on $T^*G$. One of them is conjugation by elements of $G$. For 
$g\in G$ we consider $\Conj_g:G\to G$, given by $\Conj_g(a)=g.a.g\i$, 
and the induced mapping $T^*\Conj_g = (T\Conj_{g\i})^*:T^*G\to T^*G$. 
Then we have
$$\align
(\ze^l,\pi)T^*(\Conj_g)(\xi_a) &= 
     ((T\mu^{gag\i})^*(T\mu^g.T\mu_{g\i})^*.\xi_a, g.a.g\i)\\
&= (\Ad(g\i)^* (T\mu^a)^*.\xi_a, g.a.g\i) = 
     (\Ad^*(g)\ze^l(\xi_a),g.a.g\i)\\
&= (0,g).(\ze^l(\xi_a),a).(0,g)\i.
\endalign$$
Obviously $\ze^l\o(\ze^l,\pi)\i = \pr_1:\g^*\x G\to \g^*$, and from 
$$\ze^r|T_g^*G = \Ad^*(g\i).\ze^l|T_g^*G = \Ad(g)^*.\ze^l|T_g^*G $$ 
we get that 
$$\ze^r_i(\ze^l,\pi)\i(\al,g) = \langle \al, \Ad(g).X_i\rangle.$$
Since $\pi\o(\ze^l,\pi)\i = \pr_2:\g^*\x G \to G$ we get
for both Maurer-Cartan forms $\ka^l, \ka^r$ that 
$$\pi^*\ka^{l,r}\o T(\ze^l,\pi)\i = \ka^{l,r}\o T\pi \o T(\ze^l,\pi)\i = 
     \ka^{l,r} \o T\pr_2.$$
So using \nmb!{2.3}, we get for the symplectic form
$$\align
\om^r:&= ((\ze^l,\pi)\i)^*\om 
= d\;((\ze^l,\pi)\i)^*\langle \ze^l, \pi^*\ka^r\rangle
= d\langle \pr_1, \pr_2^*\ka^r\rangle \\
&= \langle d\pr_1, \pr_2^*\ka^r\rangle 
     + \langle \pr_1, \pr_2^*d\ka^r\rangle\\
&= \langle d\pr_1, \pr_2^*\ka^r\rangle 
     + \langle \pr_1, \tfrac12\pr_2^*[\ka^r,\ka^r]_\g\rangle.\tag1
\endalign$$
The other expressions of \nmb!{2.2} for the Liouville form $\Th$ 
quickly lead again to \thetag1.

We look for an explicit expression of the Poisson structure 
$\La^r := ((\ze^l,\pi)\i)^*\La\in\Ga(\La^2T^*(\g^*\x G))$. For that we fix 
a basis $X_1,\dots, X_n$ of the Lie algebra $\g$ with structure 
constants $[X_i,X_j]=\sum_kc^k_{ij}X_k$ as in \nmb!{2.4}. We consider the  
dual basis $\xi_1,\dots,\xi_n$ in $\g^*$ with coordinate functions 
$x^1,\dots,x^n:\g^*\to \Bbb R$, so that $Id_{\g^*}=\sum x^i\xi_i$.
In these coordinates the symplectic structure is given by 
$$
\om^r= ((\ze^l,\pi)\i)^*\om = \sum_i dx^i\wedge \ka^r_i + \tfrac12 
\sum_{ijk} x^ic^i_{jk}\ka^r_j\wedge \ka^r_k\tag2
$$
By computing $\ddt (\ze^l,\pi)T^*(\Fl^{L_i}_t)(\xi_g)$ one easily 
checks that the vector field $L_i^*\in \X(T^*G)$ is 
$(\ze^l,\pi)$-related to $0\x L_i\in \X(\g^*\x G)$; for the 
prolongations of the right invariant vector fields one obtains
that $R^*_i$ is $(\ze^l,\pi)$-related to 
$\ad(-X_i)^*\x R_i = \ad^*(X_i)\x R_i\in \X(\g^*\x G)$.
Note that 
$\ad(X_i)^*(\al)=\al\o \ad(X_i) = \sum_j x^j(\al).(\xi_j\o \ad(X_i))
= \sum_{jk} x^j(\al).c^j_{ik}\xi_k$, so that the vector field 
$\ad(X_i)^*\in\X(\g^*)$ is given by 
$\sum_{jk} x^jc^j_{ik}\tfrac\partial{\partial x^k}$. 
Thus we get 
$i(\ad(-X_i)^*\x R_i)\om^r = - dx^i$ and
$i\left(-\tfrac\partial{\partial x^i}\x 0\right)\om^r = -\ka^r_i$.
The Poisson structure is given by
$$\align
\La^r &= ((\ze^l,\pi)\i)^*\La \\
&= \sum_i \ad(X_i)^*\wedge \tfrac\partial{\partial x^i} -
     \sum_i R_i\wedge \tfrac\partial{\partial x^i}
+ \tfrac12 \sum_{ijk} x^ic^i_{jk}\tfrac\partial{\partial x^j}
     \wedge \tfrac\partial{\partial x^k}\\
&=    \sum_{ijk} x^jc^j_{ik}\tfrac\partial{\partial x^k}\wedge 
     \tfrac\partial{\partial x^i}
     - \sum_i R_i\wedge \tfrac\partial{\partial x^i}
+ \tfrac12 \sum_{ijk} x^ic^i_{jk}\tfrac\partial{\partial x^j}
     \wedge \tfrac\partial{\partial x^k}\\
&=      - \sum_i R_i\wedge \tfrac\partial{\partial x^i}
- \tfrac12 \sum_{ijk} x^ic^i_{jk}\tfrac\partial{\partial x^j}
     \wedge \tfrac\partial{\partial x^k}\tag 3\\
\endalign$$
The Poisson bracket of two functions $f,g\in C^\infty(\g^*\x G)$ is 
then given by 
$$-\{f,g\} = \sum_i (0\x R_i)(f). \tfrac{\partial g}{\partial x^i}
- \sum_i (0\x R_i)(g). \tfrac{\partial f}{\partial x^i}
+ \sum_{ijk} x^ic^i_{jk}\tfrac{\partial f}{\partial x^j}
     \tfrac{\partial g}{\partial x^k}\tag4
$$

\head\totoc\nmb0{3}. Generalizing momentum mappings \endhead

\subhead\nmb.{3.1}. A more general construction \endsubhead
For a Lie group $G$ with Lie algebra $\g^*$ we generalize the formulae
obtained in section \nmb!{2}.
Instead of the left or right momentum  
let us start with an arbitrary smooth mapping $f:T^*G\to \g^*$ and 
in view of \nmb!{2.3} let us consider the 1-form
$$\Th_f = \langle f,\pi^*\ka^l\rangle \in \Om^1(T^*G).$$
The 2-form $\om_f = d\Th_f$ is closed but it may be degenerate. We 
try reduction by the kernel $\ker\om_f$ of $\om_f$ in order to get the 
symplectic form $\tilde\om_f$ on $T^*G/\ker\om_f$. This quotient 
space is a smooth manifold only if we restrict ourselves to such open 
subsets of $T^*G$ on which $\ker\om_f$ is a smooth distribution of 
constant rank.

The most important case is when $f$ is invariant with respect to 
the left action of $G$ on $T^*G$. Then $f$ may be viewed as a 
function in the coordinates $\ze^r_i$ of the right momentum $\ze^r$, so 
$f = f(\ze^r_1,\dots,\ze^r_n)$, and as in \nmb!{2.4} we get
$$\align
f &= \sum_{i=1}^n f_i.\xi_i,\\ 
\Th_f &= \sum _if_i.\pi^*\ka^l_i = \sum_i\tilde f_i.\pi^*\ka^r_i, \\
\om_f &= \sum_i\left( df_i\wedge \pi^*\ka^l_i + 
     \frac12\sum_{j,k}c^i_{jk} f_i.\ka^l_j\wedge \ka^l_i\right),\\
-i_{R_j^*}\om_f &= d\tilde f_j.
\endalign$$ 
The kernel of $\om_f$ depends on the $d f_i$'s and also on the 
structure constants, and it is difficult to describe in this 
generality.

\subhead\nmb.{3.2}. Coadjoint orbits \endsubhead
We show here that each orbit of the coadjoint action with its 
symplectic structure can be obtained by the method in \nmb!{3.1}.
For $\xi\in \g^*$ let $G_\xi$ be the isotropy group with Lie algebra 
$\g_\xi\subset \g$. We let $f=\xi$, the constant mapping 
$T^*G\to \g^*$. Then 
$$\align
\Th_\xi &= \langle \xi,\pi^*\ka^l\rangle = 
     \pi^*\langle \xi,\ka^l\rangle = \pi^*\ka^l_\xi,\\
\om_\xi &= \pi^*\langle \xi,d\ka^l\rangle = 
     -\frac12\langle \xi,[\pi^*\ka^l,\pi^*\ka^l]^\wedge\rangle \\
&= \frac12 \langle \ad(\pi^*\ka^l)^*\xi,\pi^*\ka^l \rangle, 
\endalign$$
so the kernel of $\om_\xi$ is spanned by the vertical vector fields 
and by the $L^*_X$ for $X\in \g_\xi$. Thus $T^*G/\ker\om_\xi = 
G/G_\xi$, and the last formula for $\om_\xi$ shows that the induced 
symplectic structure on $G/G_\xi$ corresponds to the canonical 
symplectic structure on the coadjoint orbit through $\xi$.

\subhead\nmb.{3.3}. Example \endsubhead
The reduction from \nmb!{3.2} is used in mechanics with respect to 
gauge invariant Lagrangians. Consider for instance the Lagrangian of 
the free spinning particle \cit!{3}, \cit!{5}
$$L=\frac12 \sum_{k=1}^3 m\dot q^2_k + \la i \tr(\si_3.g\i.\dot g)$$
where $\si_3$ is the Pauli matrix and $g\in SU(2)$, so that the 
configuration space is $\Bbb R^3\x SU(2)$.

The 1-form $\Th_L$ of this Lagrangian can be written in the form
$$\Th_L= \sum_{k=1}^3 p_k\;dq_k + \la i \tr(\si_3.g\i.dg),$$
where $p_k=m\dot q_k$ and where $SU(2)$ is just the matrix group.
Using the standard basis $X_k=\frac12 i\si_k$ in $\frak s\frak u(2)$ 
and the dual basis $\xi_k$ we may write in our previous notation
$$\align
\Th_L &= \sum_{k=1}^3 p_k\;dq_k + \la \pi^* \ka^l_3 ,\\
\om_L &= d\Th_L = \sum_{k=1}^3 dp_k\wedge dq_k 
     + \la \pi^* (\ka^l_1\wedge \ka^l_2),
\endalign$$
and the kernel of $\om_L$ is spanned by the vertical vector fields on 
$T^*SU(2)$ and by $L^*_3$. The reduced space is therefore 
$$T^*\Bbb R^3\x SU(2)/SU(1) \cong T^*\Bbb R^3 \x S^2,$$
where $SU(1)$ is the subgroup in $SU(2)$ generated by $L_3$. This is 
exactly the phase space for the free spinning particle.

\subhead\nmb.{3.4}. Reduction to a Cartan subalgebra \endsubhead
Let now $G$ be a compact Lie group with Lie algebra $\g$, and choose 
a Cartan subalgebra $\h\subset \g$. Let $\De$ be the system of roots, 
choose a positive root system $\De_+$ and denote by 
$\Pi=\{\al_1,\dots,\al_k\}$ the associated simple roots.
Let 
$$\{E_\al, E_{-\al}: \al\in\De_+\}\cup \{H_1,\dots,H_k\}$$
be the corresponding Cartan basis in $\g$ and let 
$$\{\xi_\al, \xi_{-\al}: \al\in\De_+\}\cup \{\xi_1,\dots,\xi_k\}$$
be the dual basis in $\g^*$. Let us expand the right momentum in it,
$$\ze^r = \sum_{j=1}^k\ze^r_j\xi_j + \sum_{\al\in\De}\ze^r_\al\xi_\al,$$
and let us consider 
$$f:= \sum_{j=1}^k\ze^r_j\xi_j.$$
The corresponding 2-form $\om_f$ is then given by 
$$\align
\om_f &= \sum_{j=1}^k d\ze^r_j\wedge \pi^*\ka^l_j + 
     \sum_{j=1}^k \ze^r_j\wedge \pi^*d\ka^l_j \\
&= \sum_{j=1}^k d\ze^r_j\wedge \pi^*\ka^l_j - 
     \frac12\sum_{\al\in\De_+}
     \frac{2B(f,\al)}{B( \al,\al)} 
     \pi^*\ka^l_\al \wedge \pi^* \ka^l_{-\al},
\endalign$$
where $B( \quad,\quad)$ is the Cartan Killing form on 
$\g^*$. Since $\bigcap_{j=1}^k \ker d\ze^r_j \subset \ker \om_f$ we can 
reduce $\om_f$ to a 2-form $\tilde\om_f$ on the respective quotient 
manifold which can be identified with $G\x\h^*$:
$$\tilde\om_f = \sum_{j=1}^k d\ze^r_j\wedge \ka^l_j - 
     \sum_{\al\in\De_+} \frac{B(f,\al)}{B( \al,\al)} 
     \ka^l_\al \wedge \ka^l_{-\al},$$
where $\ze^r_j$ is regarded as a coordinate just in $\h^*$. The 
2-form $\tilde\om_f$ is now non-degenerate and hence symplectic, on 
$G\x C$ for each open Weyl-chamber $C$ in $\h^*$. We also have the 
following correspondence (by dual bases)
$$
L_j \leftrightarrow d\ze^r_j, 
\quad L_\al \leftrightarrow \frac{B(f,\al)}{B(\al,\al)}\ka^l_{-\al}
     \text{ for }\al\in\De,
\quad Z^r_j \leftrightarrow -\ka^l_j,
$$
so the associated Poisson structure may be written as
$$\La = \sum_{j=1}^k Z^r_j\wedge L_j - 
     \sum_{\al\in\De_+}
     \frac{B(\al,\al)}{B(f,\al)} 
     L_\al \wedge L_{-\al}.$$
For $G=SU(2)$ for example, we get the Poisson structure
$\La=\frac\partial{\partial p}\wedge L_3 - \frac1p L_1\wedge L_2$
on $\Bbb R^*\x SU(2)$, where $L_j=i \si_j$, $j=1,\dots,3$ is the 
standard basis of $\frak s\frak u(2)$.

We observe finally, that by fixing a value $\xi$ of $f$ we obtain a 
submanifold $\Si\subset G\x \h^*$ on which $\tilde\om_f$ has kernel 
generated by all left invariant vector fields on $G$ corresponding to 
the isotropy Lie algebra $\g_\xi$, and the reduction gives again the 
canonical symplectic structure on the orbit $G/G_\xi$.

A similar formula as the last one is contained in \cit!{1}, without 
geometric interpretation, in the context of quantum groups.

\head\totoc\nmb0{4}. Symplectic structures on cotangent bundles of 
principal bundles \endhead

\subhead\nmb.{4.1}. Symplectic forms on $T^*P$ 
\endsubhead
Let $p:P\to M^n$ be a principal $G$-bundle, i\.e\. there is a free 
right action $r:P\x G\to P$ of $G$ on $P$ and $M$ is the orbit space, 
and we suppose that $p$ is a locally trivial fiber bundle. For 
$X\in \g$, the Lie algebra of $G$, we denote by $\ze_X\in\X(P)$ the 
fundamental vector field of the principal right action, 
$\ze_X(u)=\ddt u.\exp(tX)$. Then $\ze:\g\to \X(P)$ is an injective 
homomorphism of Lie algebras and it gives us an isomorphism 
$\ze_u:\g\to V_uP$ onto the vertical space 
$V_uP:= \ker(T_up:T_uP\to T_{p(u)}M) = T_u(u.G)$ for any $u\in P$. 
The inverse isomorphisms $\ka^v_u= \ze_u\i:V_uP\to \g$ form a mapping 
$\ka^v: VP\to \g$ which we call the vertical parallelism. Note for 
further use that from $T(r^g).\ze_X(u)=\ze_{\Ad(g\i)X}(u.g)$ we have 
$$\Ad(g)\o \ka^v_{u.g} = \ka^v_u\o T(r^{g\i}).\tag0$$

Let now $\ga\in\Om^1(P;\g)$ be a principal connection in $P$, i\.e\. 
a $G$-equivariant $\g$-valued 1-form on $P$ that prolongs the 
vertical parallelism $\ka^v$. So we have:
\roster
\item $(r^g)^*\ga = \Ad(g\i).\ga$ for each $g\in G$, where 
       $r^g:P\to P$ is the right action by $g$.
\item $\ga|VP = \ka^v$.
\endroster

If $X_1,\dots,X_k$ is a basis of the Lie algebra $\g$, we can expand 
$\ga$ by 
$$\ga = \sum_{j=1}^k \ga^j.X_j.\tag3$$
Let $\pi=\pi_P:T^*P\to P$ be the canonical projection and consider 
the pullback $\pi^*\ga\in \Om^1(T^*P;\g)$. Since $\pi:T^*P\to P$ is 
$G$-equivariant, the pullback $\pi^*\ga$ is a $G$-equivariant 
$\g$-valued 1-form on $T^*P$, so that 
$(T^*r^g)^*\pi^*\ga = ((Tr^{g\i})^*)^*\pi^*\ga = \pi^*(r^g)^*\ga = 
\Ad(g\i)\pi^*\ga$.

For any smooth mapping $f:T^*P\to \g^*$, in coordinates
$f=\sum_{j=1}^k f_j.\xi_j$, where $\xi_1,\dots,\xi_k$ is the basis of 
$\g^*$ dual to $(X_i)$, we may then consider the 1-form
$$
\Th_f= \Th_{f,\ga} = \langle f,\pi^*\ga\rangle_\g 
     = \sum f_j.\ga^j\in \Om^1(T^*P).
\tag4$$

We consider now $\Th_f$ as a generalized Liouville form on $T^*P$, 
corresponding to the moment mapping $f:T^*P\to \g^*$, and we consider 
the 2-form 
$$\align
\om_f:&= d\Th_f = \langle df,\pi^*\ga\rangle^\wedge 
     + \langle f,\pi^*d\ga\rangle \tag5\\
&= \langle df,\pi^*\ga\rangle^\wedge 
     + \langle f, \pi^*(\Om -\tfrac12[\ga,\ga]^\wedge)\rangle\\
&= \sum_{j=1}^k  df_j\wedge \pi^*\ga^j
     + \sum_{j=1}^k f_j.\pi^*\Om^j
     -\tfrac12\sum_{j,k,l} c^j_{kl}.f_j.\pi^*\ga^k\wedge \ga^l,
\endalign$$
where $\Om=d\ga+\frac12[\ga,\ga]^\wedge = \sum_j \Om^j X_j$ 
is the curvature form of the principal connection $\ga$.
This 2-form is closed but no longer non degenerate and we try to 
reduce it to a symplectic form on $P/\ker\om_f$. This cannot be done 
in general and we will discuss now some particular choices of $f$.

\subhead\nmb.{4.3}. The canonical momentum \endsubhead
There is a canonical choice of the momentum $f:T^*P\to \g^*$. 
Namely, 
the action of $G$ on $T^*P$ is Hamiltonian with respect to the 
canonical symplectic structure $\om_P = d\Th_P$ on $T^*P$. The 
associated canonical momentum mapping is given by
$$\gather
f_{\text{can}}: T^*P \to \g^*,\\
\langle f_{\text{can}}(\ph), X\rangle = 
     \langle \ph,\ze_X\rangle,\quad \ph\in T^*P, X\in \g,
\endgather$$
and the associated mapping $(\ka^v)^*_{f_{\text{can}}}:T^*P\to V^*P$ is the 
fiberwise adjoint of the inclusion $VP\to TP$. So this moment mapping 
is invariant under all gauge transformations, by \nmb!{4.2}.
For the associated Liouville form and its derivative we have
$$\Th^v=\Th_{f_{\text{can}}} = \langle f_{\text{can}},\pi^*\ga \rangle,$$
which we will call the \idx{\it vertical Liouville form}.

\subhead\nmb.{4.4}. The canonical momentum in coordinates
\endsubhead
We continue our investigation in a principal bundle chart now, so we 
assume that $P=\Bbb R^n\x G$ and $T^*P=T^*\Bbb R^n\x T^*G$. We will 
use coordinates $q_i$ in $\Bbb R^n$ and $(q_i,p_i)$ in $T^*\Bbb R^n$.
We use again a basis $X_i$ of the Lie algebra $\g$ with dual basis
$\xi_i$ in $\g^*$. Then $f_{\text{can}}:T^*P\to \g^*$ is of the form 
$f_{\text{can}}(\al',\al'')=\ze^r(\al'')$, where $\al'\in T^*\Bbb R^n$, 
$\al''\in T^*G$, and where $\ze^r:T^*G\to \g^*$ is the right momentum 
mapping from \nmb!{2.2}. 

The principal connection $\ga\in\Om^1(P;\g)$ 
with coordinate expression $\ga=\sum_j\ga^j X_j$ may then be written 
in terms of the vector potential $A$ as 
$$\gather
\ga_{(q,g)} = \ka^l_g + \Ad(g\i)A_q,\\
A=\sum_jA^jX_j = \sum_{i,j}A^j_i dq_i\otimes X_j\in\Om^1(\Bbb R^n;\g).
\endgather$$
The curvature form $\Om$ of $\ga$ is then expressed in terms of the 
curvature $F$ by 
$$\gather
\Om = d\ga +\tfrac12[\ga,\ga]^\wedge,\qquad
\Om_{(q,g)} = \Ad(g\i)F_q,\\
F=dA+\tfrac12[A,A]^\wedge = \sum_i F^iX_i = 
\sum_{i,j,k}F^i_{jk}\,dq_j\wedge dq_k\otimes X_i\in\Om^2(\Bbb R^n,\g).
\endgather$$ 
 From \nmb!{4.3} we have for the vertical Liouville form 
$$\align
\Th^v &= \langle f_{\text{can}},\pi^*\ga \rangle = \sum_i \ze^r_i.\pi^*\ga^i,\\
\Th^v_{(\al'_q,\al''_g)} &= \langle \ze^r(\al''_g),\ka^l_g T\pi_G + 
     \Ad(g\i)A_qT\pi_{\Bbb R^n} \rangle\\
&= \langle \Ad^*(g)\ze^r(\al''_g),\Ad(g)\ka^l_g T\pi_G + 
     A_qT\pi_{\Bbb R^n} \rangle\\
&= \langle \ze^l(\al''_g),\ka^r_g T\pi_G + 
     A_qT\pi_{\Bbb R^n} \rangle,\\
\Th^v &= \Th_{T^*G} + \langle \ze^l,\pi^* A \rangle,\tag1
\endalign$$
since we have $\ka^l_g = \Ad(g\i)\ka^r_g$ and 
$\Ad(g\i)^*\ze^r=\ze^l$, see \nmb!{2.2}.  The 1-form $\Th^v$ in 
\thetag1 was 
considered in \cit!{3}, \cit!{8}, \cit!{10}, and \cit!{11}.
Let us call
$\et=\pi^*\ka^r + \pi^* A = \sum_j\et_j\otimes X_j \in \Om^1(T^*P;\g)$, 
then using this trick once more we get the exterior 
derivative of $\Th^v$ as 
$$\align
\om^v &= d\Th^v =  d \langle\ze^r,\pi^*\ga\rangle \tag2\\
&= \langle d\ze^r,\pi^*\ga\rangle^\wedge  
     - \langle\ze^r,\pi^*[\ga,\ga]^\wedge \rangle  
     + \langle\ze^r,\pi^*\Om\rangle\\ 
&= \langle d\ze^r,\pi^*\ga\rangle^\wedge  
     - \langle\ze^l,[\et,\et]^\wedge \rangle  
     + \langle\ze^l,\pi^*F\rangle\\ 
&= \sum_i d\ze^r_i\wedge \pi^*\ga^i 
     - \frac12\sum_{ijk}\ze^l_i\, c^i_{jk}\, 
     \et_j\wedge \et_k + \sum_i\ze^l_iF^i_{jk}\, 
     dq_j\wedge dq_k. 
\endalign$$
Equation \thetag1 shows that $\Th^v$ is invariant with respect to the 
right action of $G$ on $T^*P$ induced from the principal right action 
on $T^*P$. For the infinitesimal generators $L^*_j$ of this action 
(see \nmb!{2.4}) we have then for the Lie derivative 
$$\L_{L^*_j}\Th^v = 0.\tag3$$

We add now to the vertical Liouville form a horizontal one, for which 
we choose the pullback of the canonical Liouville form 
$\Th_{\Bbb R^n}\in\Om^1(T^*\Bbb R^n)$ (We could choose a more general 1-form 
here): 
$$\Th_\ga = \Th^v + \operatorname{pr}^*\Th_{\Bbb R^n} 
     = \Th^v + \sum_i p_idq_i, \quad \om_\ga = d\Th_\ga,$$
where $(q_i,p_j)$ are the standard coordinates on $T^*\Bbb R^n$. This 
depends on the choice of the trivialization since we have no 
canonical projection $T^*P\to T^*M$.
The 2-form $\om_\ga$ is a symplectic form on 
$T^*P\cong T^*\Bbb R^n\x T^*G$, and since $\L_{L^*_j}\Th_\ga=0$, we 
get $-i_{L^*_j}\om_\ga =  d\ze^r_j$. So the isomorphism between 
1-forms and vector fields on $T^*P$ induced by $\om_\ga$ gives us the 
correspondence $L_j^*\leftrightarrow  d\ze^r_j$. Similarly as in 
\nmb!{2.4} we find vector fields $\partial_{\ze^r_j}\in\X(T^*P)$ such 
that $\partial_{\ze^r_j}\leftrightarrow -\pi^*\ga^j$. Equation \thetag1 
shows that we can find $\partial_{\ze^l_j}\leftrightarrow -\et_j$. 
It is easy to see that $\partial_{p_j}\leftrightarrow - dq_j$.
 From equation \thetag2 and the fact that
$\langle \partial_{q_j},\et_k\rangle = 
     \langle \partial_{q_j},A^k\rangle = A^k_j$
we see that 
$$\partial_{q_j}\longleftrightarrow 
     \sum_i b^i_j d\ze^r_i - 
     \sum_{i,k,s}\ze^l_ic^i_{k,s}A^s_j\et_k + 2\sum_{i,k}\ze^l_i 
     F^i_{kj} dq_k + dp_j,$$
where $b^i_j=\langle  \ga_i,\partial_{q_j}\rangle$. 
So if we put 
$$\widetilde{\partial_{q_j}} := \partial_{q_j} - 
     \sum_i b^i_j L^*_i - 
     \sum_{i,k,s}\ze^l_ic^i_{ks}A^s_j\partial_{\ze^l_k} 
     + 2\sum_{i,k}\ze^l_i F^i_{kj} \partial_{p_k},$$
we get $\widetilde{\partial_{q_j}}\leftrightarrow dp_j$. Therefore we 
can write the Poisson structure $\La_\ga$ corresponding to $\om_\ga$ 
in the form
$$\align
\La_\ga &=  \sum_i \partial_{\ze^r_i}\wedge L^*_i - 
     \frac1{2}\sum_{i,j,k}\ze^l_ic^i_{jk}
     \partial_{\ze^l_j}\wedge \partial_{\ze^l_k} 
     + \sum_{i,j,k}\ze^l_i F^i_{jk} 
     \partial_{p_j}\wedge \partial_{p_k}\tag4\\
&\quad + \sum_j\partial_{p_j}\wedge 
     \left(\partial_{q_j} - \sum_i b^i_j L^*_i - 
     \sum_{i,k,s}\ze^l_ic^i_{ks}A^s_j\partial_{\ze^l_k} 
     + 2\sum_{i,k}\ze^l_i F^i_{kj} \partial_{p_k}\right)\\
&=  \sum_i \partial_{\ze^r_i}\wedge L^*_i - 
     \sum_{i,j} b^i_j \partial_{p_j}\wedge L^*_i
     - \frac1{2}\sum_{i,j,k}\ze^l_ic^i_{jk} 
     \partial_{\ze^l_j}\wedge \partial_{\ze^l_k}\\
&\quad - \sum_{i,j,k}\ze^l_i F^i_{jk} 
     \partial_{p_j}\wedge \partial_{p_k}
     - \sum_{i,j,k,s}\ze^l_ic^i_{ks}A^s_j 
     \partial_{p_j}\wedge\partial_{\ze^l_k} 
     + \sum_{j}\partial_{p_j}\wedge\partial_{q_j}.
\endalign$$
Since $\La_\ga$ is invariant with respect to the $G$-action we can 
reduce $\La_\ga$ to $\tilde\La_\ga$ on 
$T^*P/G\cong T^*\Bbb R^n\x(T^*G/G)\cong T^*\Bbb R^n\x \g^*$ 
by putting $L_j^*=0$ and considering $I_k:=\ze^l_k$ as coordinate in 
$\g^*$:
$$\align
\tilde\La_\ga &= 
     - \frac1{2}\sum_{i,j,k}I_ic^i_{jk}
     \partial_{I_j}\wedge \partial_{I_k} 
     - \sum_{i,j,k}I_i F^i_{jk} 
     \partial_{p_j}\wedge \partial_{p_k}\tag5\\
&\quad -  \sum_{i,j,k,s}I_ic^i_{ks}A^s_j \partial_{p_j}\wedge\partial_{I_k} 
     + \sum_{j}\partial_{p_j}\wedge\partial_{q_j}.
\endalign$$
It should be noticed that this bivector field is degenerate, 
therefore it does not have a symplectic description on $T^*P/G$. In 
this respect our result is quite different from those available in 
the literature. Our description deals with several particles instead 
of a given one with fixed isospin or color. That latter is obtained 
by fixing a value for the Casimir functions of our brackets. In the 
generic case we find a symplectic level set diffeomorphic to 
$T^*\Bbb R^n \x \Ad^*(G)c$, where $\Ad^*(G)c$ is the coadjoint 
orbit through $c$ in $\g^*$.
The Poisson structure $\tilde\La_\ga$ 
depends on the choice of the connection $\ga$ 
and its curvature $\Om$, so in the local trivialization on 
the vector potential $A$ and its curvature $F$. If for example 
instead of equation \thetag1 we consider 
$\Th^v = \Th_{T^*G} + \la \langle \ze^l,\pi^* A \rangle$ 
(where $\la$ could be absorbed into 
the choice of $A$, so it is not more general) then in the expression 
\thetag4 for $\La_\ga$ the constant $\la$ would appear as a coupling 
constant. The Poisson bracket associated to \thetag5 appears already in 
\cit!{9}. 

Let us find the Hamiltonian vector field $\Ga$ corresponding to the 
free Hamiltonian
$$H=\frac12\sum_j p_j^2.$$
It turns out to be 
$$\Ga = 2\sum_{i,k,j}I_i  F^i_{kj} p_j\partial_{p_k} - 
\sum_{i,j,k,s}I_ic^i_{ks}A^s_jp_j\partial_{I_k} + 
\sum_jp_j\partial_{q_j},$$
which describes the motion of a Yang-Mills particle which carries a 
`charge' given by the spin-like variable $I_k$. 
In the free case, $A=0$, we would have 
$$\align
\tilde\La_\ga &= 
     - \frac1{2}\sum_{i,j,k}I_ic^i_{jk}
     \partial_{I_j}\wedge \partial_{I_k} 
     + \sum_{j}\partial_{p_j}\wedge\partial_{q_j},\\
H &=\frac12\sum_j p_j^2.
\endalign$$
It is interesting to 
notice that the presence of the vector potential changes the Poisson 
bracket without changing the Hamiltonian.

When $G$ is the group $U(1)$ the Yang-Mills field (curvature $F$) 
reduces to electromagnetism and the vector field $\Ga$ contains the 
standard Lorentz force expression:
$$\Ga = 2\sum_{k,j}e F_{kj} p_j\partial_{p_k} + \sum_jp_j\partial_{q_j},$$
where $e$ is the electric charge and the equation of motion in the 
internal variables reduces to $\frac{\partial e}{\partial t}=0$.

\subhead\nmb.{4.2}. Behavior under gauge transformations \endsubhead
Let $\operatorname{Gau}(P)$ denote the group of all gauge 
transformations of the principal bundle $p:P\to M$. So elements of 
$\operatorname{Gau}(P)$ are diffeomorphisms $\ph:P\to P$ which 
respect fibers $(p\o \ph=p)$ and which commute with 
the principal right action $r$ of $G$ on $P$.

If $\ga\in\Om^1(P;\g)$ is a principal connection in $P$, let 
$\ze_\ga:TP\to VP$ be the $G$-equivariant projection onto the 
vertical bundle induced by $\ga$, i\.e\. 
$\ze_\ga(\xi_u)=\ze_{\ga(\xi_u)}(u)$, where $\xi_u\in T_uP$ and 
$\ze:\g\to \X(P)$ is the fundamental vector field mapping of the 
principal right action. A gauge transformation 
$\ph$ act naturally on a connection $\ga$:
$$
\ph^*(\ze_\ga):= T\ph\i\o \ze_\ga\o T\ph = \ze_{\ph^*\ga}
     = \ze_{\ga\o T\ph},
\tag 1$$
where the second equation follows from properties of $\ze$.
We may describe a gauge transformation $\ph\in\operatorname{Gau}(P)$ 
either as a section $s_\ph\in\Ga(P\x_GG)$ of the 
associated bundle $P\x_GG$, where the 
structure group $G$ acts on the fiber $G$ by conjugation, or 
equivalently as a $G$-equivariant mapping 
$f_\ph\in C^\infty(P,(G,\operatorname{Conj}))^G$. 
Then $\ph(u)=u.f_\ph(u)$, and the isomorphism 
$\Ga(P[G,\operatorname{Conj}])\cong 
     C^\infty(P,(G,\operatorname{Conj}))^G$ 
is standard, see e\.g\. \cit!{4},~10.14.
Then we get for the action of $\ph$ on a connection $\ga$
$$
\ph^*\ga= \ga\o T\ph = \Ad(f_\ph\i)\ga + f_\ph^*\ka^l,
\tag2$$ 
where $f_\ph^*\ka^l$ is the the pullback of the left Maurer-Cartan 
form on $G$. For the curvature forms one has 
$$
\Om_{\ph^*\ga}=\Om_{\Ad(f_\ph\i)\ga + f_\ph^*\ka^l} = 
     \ph^*\Om_\ga = \Ad(f_\ph\i)\Om_\ga.
$$
For a smooth mapping $f:T^*P\to\g^*$ the action of a gauge 
transformation $\ph\in \operatorname{Gau}(P)$ on the 1-form
$\Th_{f,\ga}=\langle f,\pi^*\ga\rangle$ from \nmb!{4.1},~\thetag4 
turns out as
$$ (T^*\ph)^*\Th_{f,\ga} 
     = \langle f\o T^*\ph,\pi^*(\Ad(f_\ph\i)\ga + f_\ph^*\ka^l)\rangle 
     = \Th_{f\o T^*\ph,\ph^*\ga}.
\tag3$$
Using the vertical parallelism $\ka^l$ we can associate to 
$f:T^*P\to \g^*$ the fiber respecting smooth mapping 
$(\ka^v)^*_f:T^*P\to V^*P$ which is given by 
$((\ka^v)^*_f)|T^*_uP = 
(\ka^v_u)^*\o (f|T^*_uP):T^*_uP\to\g^*\to V^*_uP$. Note that 
$(\ka^v)^*_f$ is 
$G$-equivariant, i\.e\. $(k^v)^*_f\o (T^*r^g) = (T^*r^g)\o (k^v)^*_f$ 
for all $g\in G$, if and only if  
$f\o (T^*r^g) = \Ad^*(g\i)\o f=\Ad(g)^*\o f$. 
We can then express $\Th_{f,\ga}$ also by
$$
\Th_{f,\ga} = \langle f,\pi^*\ga\rangle_\g 
     = \langle (k^v)^*_f,\ze_\ga\o T\pi \rangle_{VP}
\tag4$$
from which is follows very easily that $\Th_{f,\ga}$ is $G$-invariant 
if and only if the mapping $(\ka^v)^*_f:T^*P\to V^*P$ is $G$-equivariant, i\.e\. 
$f\o (T^*r^g) = \Ad^*(g\i)\o f=\Ad(g)^*\o f$.

Also from $\ze_X\o\ph=T\ph\o\ze_X$ it follows in turn that 
$$\gather
(\ka^v_{\ph(u)})^*= T^*\ph\o (\ka^v_u)^*\\
(T^*\ph)^*((\ka^v)^*_f) = T^*\ph\i\o (\ka^v)^*_f\o T^*\ph 
     = (\ka^v)^*_{f\o T^*\ph},
\endgather$$
so that again all possible actions of gauge transformations on moment 
mappings $f:T^*P\to\g^*$ coincide. Note that the canonical momentum 
$f$ from \nmb!{4.3} is invariant under all gauge transformations.

Let us continue now in a local trivialization as in \nmb!{4.4}, so 
that we assume $P=\Bbb R^n\x G$. In this case a gauge transformation 
$\ph\in \operatorname{Gau}(\Bbb R^n\x G)$ is given by 
$\ph(q,g)=(q,s_\ph(q).g)$ for $s=s_\ph\in C^\infty(\Bbb R^n,G)$. If a 
connection $\ga$ is given in terms of a vector potential $A$ by 
$\ga_{(q,g)} = \ka^l_g + \Ad(g\i)A_q$ as in \nmb!{4.4}, then the 
action of the gauge transformation $\ph$ on $\ga$ is given by 
$$(\ph^*\ga)_{(q,g)} = \ka^l_g + \Ad(g\i)(s^*\ka^l + 
\Ad(s(q)\i)A_q),$$
so that $s$ acts on the vector potential and the curvature by 
$$\align
A &\mapsto s^*\ka^l + \Ad(s\i)A \\
F &\mapsto \Ad(s\i)F
\endalign$$
Now it remains the task to write down the behavior of the main 
expression in \nmb!{4.4} under gauge transformations. For that we 
have to use $\Ad(s\i).X_i=\sum_js_i^jX_j$ and to use a matrix 
representation of $s$. Note that this is not the usual way to write 
gauge transformations in a linear group $G\subset GL(N,\Bbb R)$. So 
the resulting formulas will look quite unfamiliar. But using matrix 
representations of all objects would hide the intrinsic symmetry of 
our approach and will lead to a sea of indices. 

\subhead\nmb.{4.5}. Further reduction to a Cartan subalgebra 
\endsubhead
Let us consider again the trivial principal bundle $P=\Bbb R^n\x G$ 
with compact structure group $G$, and choose a Cartan subalgebra 
$\h\subset \g$. In the notation of \nmb!{3.4}, 
instead of the momentum $T^*P\to\g^*$
pulled back from the right momentum
$$\ze^r = \sum_{j=1}^k\ze^r_j\xi_j + \sum_{\al\in\De}\ze^r_\al\xi_\al$$
on $G$, we use again the reduced momentum
$$f:= \sum_{j=1}^k\ze^r_j\xi_j.$$
For a principal connection $\ga$ on $P$ we 
consider the following Liouville form and its derivative
$$\align
\Th^\h_\ga &= \langle f,\pi^*\ga\rangle + \Th_{\Bbb R^n} \\
&= \sum_{j=1}^k\ze^r_j\ga^j + \sum_{i=1}^np_i\,dq_i, \\
\om^\h_\ga &= d\Th^\h_\ga.
\endalign$$
But now $\Th^\h_\ga$ is no longer invariant with respect to the whole 
of $G$, but it remains invariant with respect to the Cartan subgroup 
$H$ corresponding to $\h$, so
$$\L_{L^*_j}\Th^\h_\ga = 0,\quad -i_{L^*_j}\om^\h_\ga = 
     d\ze^r_j.$$
Similarly as in \nmb!{3.4} we get 
$$\align
\om^\h_\ga &= \sum_{j=1}^k d\ze^r_j\wedge \pi^*\ga^j - 
     \sum_{\al\in\De_+}\frac{B(f,\al)}{B( \al,\al)} 
     \pi^*\ga^\al \wedge \pi^* \ga^{-\al}
     + \sum_{i,j,k}\ze^r_i \Om^i_{k,j} dq_k \wedge dq_j\\ 
&\quad+ \sum_j dp_j\wedge dq_j,
\endalign$$
which after reduction can be regarded as a 2-form on 
$T^*\Bbb R^n\x G\x\h^*$. It is non degenerate on 
$T^*\Bbb R^n\x G\x C$ for any open Weyl chamber $C$.
There this symplectic form gives us the correspondence 
$$\gather
L_j \leftrightarrow d\ze^r_j, 
     \qquad L_\al \leftrightarrow 
     \frac{B(f,\al)}{B( \al,\al)} \ga^{-\al},\\
\partial_{p_i} \leftrightarrow -dq_i, 
     \qquad \partial_{\ze^r_j} \leftrightarrow -\ga^j,\\
\partial_{q_i} \leftrightarrow \sum_i b^i_j d\ze^r_i + 
     \sum_{\al\in\De}\frac{B(f,\al)}{B(\al,\al)}b^\al_j\ga^{-\al}
     +\sum_{i,k}\ze^r_i \Om^i_{k,j} dq_k + dp_j.
\endgather$$
Hence we have 
$$dp_j \leftrightarrow \partial_{q_j} 
     - \sum_ib^i_jL_i -\sum_{\al\in\De}b^\al_j L_\al
     +2\sum_{i,j,k}\ze^r_i\Om^i_{kj}\partial_{p_k},$$
where $b^\al_j=\langle  \ga^\al,\partial_{q_j}\rangle$, and we get 
the corresponding Poisson structure $\La^\h_\ga$ on 
$T^*\Bbb R^n\x G\x \h^*$ in the form 
$$\align
\La^\h_\ga &= \sum_i\partial{\ze^r_i}\wedge L_i 
     - \sum_{\al\in\De_+}\frac{B(\al,\al)}{B(f,\al)}
     L_\al\wedge L_{-\al}\\
&\quad-\sum_{i,k}\ze^r_i \Om^i_{k,j}\partial_{p_k}\wedge \partial_{p_j}
     - \sum_ib^i_j\partial_{p_j}\wedge L_i 
     -\sum_{\al\in\De}b^\al_j \partial_{p_j}\wedge L_\al
     +\sum_j \partial_{p_j}\wedge \partial_{q_j}.
\endalign$$
This Poisson structure is $L_i$-invariant, so we can do the reduction 
once more  and get the Poisson structure 
$$\align
\tilde\La^\h_\ga &= -\sum_{\al\in\De_+}
     \frac{B(\al,\al)}{B(f,\al)} L_\al\wedge L_{-\al}
     -\sum_{i,k}\ze^r_i 
     \Om^i_{k,j}\partial_{p_k}\wedge \partial_{p_j}\\
&\quad -\sum_{\al\in\De}b^\al_j \partial_{p_j}\wedge L_\al
     +\sum_j \partial_{p_j}\wedge \partial_{q_j}.
\endalign$$
on $T^*\Bbb R^n\x G/H\x \h^*$, which also gives a Poisson structure 
on $T^*\Bbb R^n\x G/H$ for any fixed value of $f$. This generalizes 
\nmb!{3.4}.

\subhead\nmb.{4.6}. Generalization to Cartan connections \endsubhead
Assume that we have now not only a connection $\ga$ but also a 
displacement form $\th$ on the principal bundle $p:P\to M$, that is a 
$G$-equivariant form $\th:TP\to V=\Bbb R^n$ with $\ker\th = VP$, 
where $\dim M=n$. Then 
$$\ka = \th +\ga: TP\to V\oplus \g$$
is a Cartan connection, i\.e\. a $G$-equivariant absolute parallelism 
on $P$. So we obtain a $(V\oplus \g)$-valued 1-form on $T^*P$:
$$\pi^*\ka = \ka\o T\pi: T(T^*P)\to TP \to V\oplus \g.$$
Now any smooth mapping 
$$f:T^*P \to V^*\oplus \g^*$$
defines a 1-form 
$$\Th_f = \langle f,\pi^*\ka\rangle\in\Om^1(T^*P)$$
which is $G$-invariant if $f$ is $G$-equivariant.

There is a canonical function $f$ as above: namely for $u\in P$ the 
Cartan connection $\ka_u:T_uP\to V\oplus \g$ is a linear isomorphism 
and 
$$f_u := (\ka_u\i)^*: T_u^*P \to V^*\oplus \g^*$$
is $G$-equivariant. 
Then $\Th_f= \Th_{P}$ is the canonical Liouville form on the 
cotangent bundle $T^*P$.

\enddocument